\documentclass[final,reqno,a4paper]{amsart}

\usepackage{graphicx}
\usepackage{nicefrac}
\usepackage{bm}
\usepackage{a4wide}
\usepackage{cite}
\usepackage{algorithm,algorithmicx,algpseudocode}
\usepackage{booktabs}
\usepackage{enumerate}

\newcommand{\dd}{\mathsf{d}}
\newcommand{\dx}{\, \dd x}

\newcommand{\NN}[1]{\left\|#1\right\|}
\newcommand{\tol}{\mathtt{tol}}
\newcommand{\Q}{\mathtt{Q}}
\newcommand{\F}{\mathsf{F}}

\DeclareMathOperator{\sech}{sech}
\newcommand{\subs}{\mathtt{subs}}
\newcommand{\subsnew}{\mathtt{subsnew}}
\newcommand{\x}[2]{\widehat{x}_{#1,#2}}
\newcommand{\xx}[1]{\widehat{\bm x}_{#1}}
\newcommand{\w}[2]{w_{#1,#2}}
\newcommand{\ww}[1]{{\bm w}_{#1}}
\newcommand{\p}{\mathtt{p}}
\newcommand{\pnew}{\mathtt{pnew}}
\newcommand\bigzero{\makebox(0,0){\text{\huge{\bf 0}}}}

\renewcommand{\sp}{\phantom{1ex}}
\newcommand{\matlab}{\textsc{Matlab}}

\newtheorem{assumption}{Assumption}

\newcommand{\myStateDouble}[1]{\State\parbox[t]{\dimexpr\linewidth-\algorithmicindent-\algorithmicindent}{#1\strut}}

\newtheorem{Remark}[equation]{Remark}

\numberwithin{equation}{section}

\numberwithin{equation}{section}

\title[Adaptive Variable Order Quadrature]{An Adaptive Variable Order Quadrature Strategy}

\author[P. Houston]{Paul Houston}
\address{School of Mathematical Sciences, University of Nottingham,
  University Park, Nottingham, NG7 2RD,
  UK}
\email{Paul.Houston@nottingham.ac.uk}

\author[T. P. Wihler]{Thomas P.~Wihler}
\address{Mathematisches Institut, Universit\"at Bern, Sidlerstrasse 5,
  CH-3012 Bern, Switzerland}
\email{wihler@math.unibe.ch}
\thanks{TW acknowledges the financial support by the Swiss National Science Foundation (SNF)}

\keywords{}

\subjclass[2010]{}

\begin{document}

\begin{abstract}
In this article we propose a new adaptive numerical quadrature procedure which includes both local subdivision of the integration domain, as well as local variation of the number of quadrature points employed on each subinterval. In this way we aim to account for local smoothness properties of the function to be integrated as effectively as possible, and thereby achieve highly accurate results in a very efficient manner. Indeed, this idea originates from so-called $hp$-version finite element methods which are known to deliver high-order convergence rates, even for nonsmooth functions. 
\end{abstract}

\subjclass[2010]{65D30,65N30}
\keywords{Adaptive quadrature, $hp$-adaptivity, adaptive Gauss quadrature.}


\maketitle


\section{Introduction}

Numerical integration methods have witnessed a tremendous development over the last few decades; see, e.g., \cite{Press:07,DahlquistBjorck:08,DavisRabinowitz:07}. In particular, adaptive quadrature rules have nowadays become an integral part of many scientific computing codes. Here, one of the first yet very successful approaches is the application of adaptive Simpson integration or the more accurate Gauss-Kronrod procedures (see, e.g., \cite{GanderGautschi:00}). The key points in the design of these methods are, first of all, to keep the number of function evaluations low, and, secondly, to divide the domain of integration in such a way that the features of the integrand function are appropriately and effectively accounted for. 

The aim of the current article is to propose a complementary adaptive quadrature approach that is quite different from previous numerical integration schemes. In fact, our work is based on exploiting ideas from $hp$-type adaptive finite element methods (FEM); cf.~\cite{HoustonSuliHPADAPT,De07,FankhauserWihlerWirz:14,MelenkAPOST01,opac-b1101124}. These schemes accommodate and combine both traditional low-order adaptive FEM and high-order (so-called spectral) methods within a single unified framework. Specifically, their goal is to generate discrete approximation spaces which allow for both adaptively refined subdomains, as well as locally varying approximation orders. In this way, the $hp$-FEM methodology is able to resolve features of an underlying unknown analytical solution in a highly efficient manner. In fact, this approach has proved to be enormously successful in the context of numerically approximating solutions of differential equations, and has been shown to exhibit high-order algebraic or exponential convergence rates even in the presence of local singularities; cf.~\cite{Schwab98,GuiBabuska86,SchotzauSchwabWihler:15}.

With this in mind, we adopt the $hp$-adaptive finite element strategy for the purpose of introducing a variable order adaptive quadrature framework. More precisely, we propose a procedure whereby the integration domain will be subdivided adaptively in combination with a local tuning of the number of quadrature points employed on each subinterval. To drive this refinement process, we employ a smoothness estimation technique from~\cite{FankhauserWihlerWirz:14,W11} (see also~\cite{HoustonSuliHPADAPT} for a related strategy), which was originally introduced in the context of $hp$-adaptive FEM. Specifically, the smoothness test makes it possible to gain local information concerning the regularity of the integrand function, and thereby, to suitably subdivide the integration domain and select an appropriate number of quadrature points for each subinterval. By means of a series of numerical experiments we demonstrate that the proposed adaptive quadrature strategy is capable of generating highly accurate approximations at a very low computational cost. The main ideas on this new approach together with a view on practical aspects will be discussed in the subsequent section.

\section{An $hp$-Type Quadrature Approach}

\subsection{General Quadrature Rules}

Typical quadrature rules for the approximation of an integral
\begin{equation}\label{eq:I}
I:=\int_{-1}^1 f(x)\dx
\end{equation}
of a continuous function~$f:\,[-1,1]\to\mathbb{R}$, take the form
\begin{equation}\label{eq:quad}
I\approx \widehat Q_p(f):=\sum_{k=1}^p \w{p}{k}f(\x{p}{k}),
\end{equation}
where~$p\ge 1$ is a (typically prescribed) integer number, and~$\{\x{p}{k}\}_{k=1}^p\subset[-1,1]$ and~$\{\w{p}{k}\}_{k=1}^p\subset(0,2]$ are appropriate quadrature points and weights, respectively. When dealing with a variable number~$p$ of quadrature points and weights, we can consider one-parameter families of quadrature rules (such as, for example, Gauss-type quadrature methods); here, for each~$p\in\mathbb{N}$, with~$p\ge p_{\min}$, where~$p_{\min}$ is a minimal number of points, there are (possibly non-hierarchical) families of quadrature points~$\xx{p}=\{\x{p}{k}\}_{k=1}^p$, and weights~$\ww{p}=\{\w{p}{k}\}_{k=1}^p$.

On an arbitrary bounded interval~$[a,b]$, $a<b$, a corresponding integration formula can be obtained, for instance, by means of a simple affine scaling
\begin{equation}\label{eq:phi}
\phi_{[a,b]}:\,[-1,1]\to[a,b],\qquad 
\widehat x\mapsto x=\phi_{[a,b]}(\widehat x)=\frac12h\widehat x+\frac12(a+b),
\end{equation}
with~$h=b-a>0$. Indeed, in this case
\[
\int_a^bf(x)\dx \approx Q_{[a,b],p}(f):=\frac{h}{2}\sum_{k=1}^p \w{p}{k}(f\circ\phi_{[a,b]})(\x{p}{k}),
\]
where~$f:\,[a,b]\to\mathbb{R}$ is again continuous. As before, for any specific family of quadrature rules, the corresponding quadrature point families~$\bm x_p$ are obtained in a straightforward way by letting $\bm x_p=\phi_{[a,b]}(\xx{p})$ (with the understanding that~$\phi_{[a,b]}$ is extended componentwise to vectors).

Furthermore, the above construction allows us to define composite quadrature rules, whereby the integral of $f$ is approximated on a collection of~$n\ge1$ disjoint (open) subintervals~$\{K_i\}_{i=1}^n$ of~$[a,b]$ with~$[a,b]=\bigcup_{i=1}^n \overline{K}_i$, i.e.,
\[
I\approx \sum_{i=1}^n Q_{K_i,p}(f|_{K_i}).
\]
In practical applications the subintervals are usually either of uniform size~$\nicefrac{(b-a)}{n}$, for sufficiently large~$n$, or alternatively, they are selected adaptively with the aim of resolving the relevant features of the given function~$f$.

\subsection{The Basic Idea: $hp$-Adaptivity}\label{sc:basic}

Adaptive quadrature rules usually generate a sequence of repeatedly bisected and possibly non-uniform subintervals~$\{K_i\}_{i=1}^n$, $n \geq 1$, of the integration domain~$[a,b]$ (i.e., each subinterval~$K_i$ may have a different length~$h_i$), with a prescribed and uniform number~$p$ of quadrature points on each subinterval. With the aim of providing highly accurate approximations with as little computational effort as possible, the novelty of the approach presented in this article is to design an adaptive quadrature procedure, which, in addition to subdividing the original interval~$[a,b]$ into appropriate subintervals, is able to adjust the number of quadrature points~$p_i$ \emph{individually} within each subinterval~$K_i$ in an effective way. We note that this idea originates from approximation theory~\cite{Scherer:80/81,DeVoreScherer:80} (see also~\cite{GuiBabuska86}), and has been applied with huge success in the context of finite element methods for the numerical approximation of differential equations. Indeed, under certain conditions, the judicious combination of subinterval refinements ($h$-refinement) and selection of local approximation orders ($p$-refinement), which results in the class of so-called $hp$-finite element methods, is able to achieve high-order algebraic or exponential rates of convergence, even for solutions with local singularities; see, e.g. \cite{Schwab98}. In an effort to automate the combined $h$- and $p$-refinement process, a number of $hp$-adaptive finite element approaches have been proposed in the literature; see, e.g, the survey article~\cite{mitchell_mcclain} and the references cited therein. In the current article, we pursue the smoothness estimation approach developed in~\cite{FankhauserWihlerWirz:14,W11} (cf.~also~\cite{HoustonSuliHPADAPT} for a related methodology), and translate the idea into the context of adaptive variable order numerical quadrature.

Starting from a subinterval~$K_i$ with~$p_i$ quadrature points, we are given a current approximation~$Q_{K_i,p_i}(f|_{K_i})$ of the subintegral
\begin{equation}\label{eq:intKi}
\int_{K_i}f(x)\dx\approx Q_{K_i,p_i}(f|_{K_i}).
\end{equation}
Then, with the aim of improving the approximate value~$Q_{K_i,p_i}(f|_{K_i})$, in the sense of an $hp$-adaptive finite element methodology in one-dimension, we propose two possible refinements of~$K_i$:
\begin{enumerate}[(i)]
\item\label{item:i} $h$-refinement: The subinterval~$K_i$ of length~$h_i$ is bisected into two subintervals~$K_i^1$ and~$K_i^2$ of equal size~$\nicefrac{h_i}{2}$, and the number~$p_i$ of quadrature points is either inherited to both subintervals or, in order to allow for derefinement with respect to the number of local quadrature points, reduced to~$p_i-1$ points. In the latter case, we obtain a potentially improved approximation
\begin{equation}\label{eq:Qh}
Q_{K_i}^{\rm h}(f)=Q_{K_i^{1},\max(1,p_i-1)}(f)+Q_{K_i^{2},\max(1,p_i-1)}(f)
\end{equation}
of~\eqref{eq:intKi}.
\item $p$-refinement: The subinterval~$K_i$ is retained, and the number~$p_i$ of quadrature points~$p_i$ is increased by~1, i.e., $p_i\gets p_i+1$. This yields an approximation
\begin{equation}\label{eq:Qp}
Q_{K_i}^{\rm p}(f)=Q_{K_i,p_i+1}(f).
\end{equation}
In case that~$p_i=p_{\max}$, where~$p_{\max}$ is a prescribed maximal number of quadrature points on each subinterval, we define
\begin{equation}\label{eq:Qp'}
Q_{K_i}^{\rm p}(f)=Q_{K_i^1,p_i}(f)+Q_{K_i^2,p_i}(f),
\end{equation}
where~$K_i^1$ and~$K_i^2$ result from subdividing~$K_i$ as in~\eqref{item:i}.
\end{enumerate}
In order to determine which of the above refinements is more appropriate for a given subinterval~$K_i$, we apply a smoothness estimation idea as outlined in the subsequent section. Once a decision between $h$- and $p$-refinement for~$K_i$ has been made, the procedure is repeated iteratively for any subintervals~$K_i$ for which~$Q_{K_i,p_i}(f|_{K_i})$ and its refined value (resulting from the chosen refinement) differ by at least a prescribed tolerance~$\mathtt{tol}>0$.

\subsection{Smoothness Estimation}\label{sc:smoothness}
The basic idea presented in the articles~\cite{FankhauserWihlerWirz:14,W11,HoustonSuliHPADAPT} is to estimate the regularity of a function to be approximated locally. Then, following along the lines of the $hp$-approximation approach, if the function is found to be smooth, according to the underlying regularity estimation test, then a $p$-refinement is performed, otherwise an $h$-refinement is employed. In~\cite{FankhauserWihlerWirz:14}, the following smoothness indicator, for a (weakly) differentiable function~$f$ on an interval~$K_j$, has been introduced (cf.~\cite[Eq.~(3)]{FankhauserWihlerWirz:14}):
\begin{equation}\label{eq:F}\tag{F}
\mathcal{F}_{K_j}[f]:=\begin{cases}\displaystyle
\frac{\NN{f}_{L^\infty(K_j)}}{h_j^{-\nicefrac12}\NN{f}_{L^2(K_j)}+\frac{1}{\sqrt2}h_j^{\nicefrac12}\NN{f'}_{L^2(K_j)}}  & \text{if }f|_{K_j}\not\equiv 0,\\[3ex]
1 & \text{if }f|_{K_j}\equiv 0.
\end{cases}
\end{equation}
The motivation behind this definition is the continuous Sobolev embedding~$W^{1,2}(K_j)\hookrightarrow L^\infty(K_j)$, which implies that
\[
\sup_{v\in H^1(K_j)}\frac{\NN{v}_{L^\infty(K_j)}}{h_j^{-\nicefrac12}\NN{v}_{L^2(K_j)}+\frac{1}{\sqrt2}h_j^{\nicefrac12}\NN{v'}_{L^2(K_j)}}\le 1;
\]
see~\cite[Proposition~1]{FankhauserWihlerWirz:14}. In particular, it follows that~$\mathcal{F}_{K_j}[f]\le 1$ in~\eqref{eq:F}; $f$ is classified as being smooth on~$K_j$ if~$\mathcal{F}_{K_j}[f]\ge\tau$, for a prescribed smoothness testing parameter~$0<\tau<1$, and nonsmooth otherwise.

To begin, we first consider the special case when~$f$ is a polynomial of degree~$p_j\ge 1$. Then, the derivative~$f^{(p_j-1)}$ of order~$p_j-1$ of $f$ is a linear polynomial, and the evaluation of the smoothness indicator~$\mathcal{F}_{K_j}\left[f^{(p_j-1)}\right]$ from~\eqref{eq:F} is simple to obtain. In fact, let us write~$f|_{K_j}$ in terms of a (finite) Legendre series, that is,
\begin{equation}\label{eq:fleg}
f|_{K_j}=\sum_{l=0}^{p_j}a_l(\widehat{L}_l\circ\phi^{-1}_{K_j}),
\end{equation}
for coefficients~$a_0,\ldots,a_{p_j}\in\mathbb{R}$. Here, $\widehat{L}_l$, $l\ge 0$, are the Legendre polynomials on~$[-1,1]$ (scaled such that~$\widehat L_l(1)=1$ for all~$l\ge 0$), and $\phi_{K_j}$ is the affine scaling of~$[-1,1]$ to~$K_j$; cf.~\eqref{eq:phi}. For~$f$ as in~\eqref{eq:fleg} it can be shown that
\begin{equation}\label{eq:Fp}
\mathcal{F}_{K_j}\left[f^{(p_j-1)}\right]=\frac{1+\xi_{p_j}}{\sqrt{1+\frac13\xi_{p_j}^2}+\sqrt2\xi_{p_j}},
\end{equation}
where~$\xi_{p_j}=(2p_j-1)\left|\nicefrac{a_{p_j}}{a_{p_j-1}}\right|$ (provided that~$a_{p_j-1}\neq 0$); see~\cite[Proposition~3]{FankhauserWihlerWirz:14}. In particular, this implies that
\begin{equation}\label{eq:range}
\frac12\approx\frac{\sqrt{3}}{\sqrt{6}+1}\le\mathcal{F}_{K_j}\left[f^{(p_j-1)}\right]\le 1;
\end{equation}
cf.~\cite[\S2.2]{FankhauserWihlerWirz:14}. 

In the context of the numerical integration rule~\eqref{eq:quad}, the above methodology can be adopted as follows: suppose we are given~$p_j\ge 2$ quadrature points and weights, $\{\x{p_j}{k}\}_{k=1}^{p_j}$ and~$\{\w{p_j}{k}\}_{k=1}^{p_j}$, respectively. Then,
\begin{equation}\label{eq:QKj}
\int_{K_j}f(x)\dx\approx Q_{K_j,p_j}(f|_{K_j})=\frac{h_j}{2}\sum_{k=1}^{p_j}\w{p_j}{k}(f\circ\phi_{K_j})(\x{p_j}{k}).
\end{equation} 
We denote the uniquely defined interpolating polynomial of~$f$ of degree~$p_j-1$ at the given quadrature points by
\[
\Pi_{K_j,p_j-1}f=\sum_{l=0}^{p_j-1}b_l(\widehat{L}_l\circ\phi^{-1}_{K_j}).
\]
Due to orthogonality of the Legendre polynomials, we note that
\[
b_l=\frac{2l+1}{h_j}\int_{K_j}\Pi_{K_j,p_j-1}f(x)(\widehat{L_l}\circ\phi^{-1}_{K_j})(x)\dx,\qquad l=0,\ldots,p_j-1.
\]
We further assume that the quadrature rule under consideration is exact for all polynomials of degree up to~$2p_j-2$. Thereby,
\begin{align*}
b_l&=\frac{2l+1}{2}\sum_{k=1}^{p_j}\w{p_j}{k}(\Pi_{K_j,p_j-1}f)\circ\phi_{K_j}(\x{p_j}{k})\widehat{L_l}(\x{p_j}{k})\\
&=\frac{2l+1}{2}\sum_{k=1}^{p_j}\w{p_j}{k}(f\circ\phi_{K_j})(\x{p_j}{k})\widehat{L_l}(\x{p_j}{k}).
\end{align*}
Consequently, we infer that
\begin{equation}\label{eq:xi}
\begin{split}
\xi_{K_j,p_j-1}:&=(2p_j-3)\left|\frac{b_{p_j-1}}{b_{p_j-2}}\right|\\
&=(2p_j-1)\frac{\sum_{k=1}^{p_j}\w{p_j}{k}(f\circ\phi_{K_j})(\x{p_j}{k})\widehat{L}_{p_j-1}(\x{p_j}{k})}{\sum_{k=1}^{p_j}\w{p_j}{k}(f\circ\phi_{K_j})(\x{p_j}{k})\widehat{L}_{p_j-2}(\x{p_j}{k})},
\end{split}
\end{equation}
and thus, in view of~\eqref{eq:Fp}, we use the quantity
\begin{equation}\label{eq:sind}
\F_{K_j,p_j}(f):=\frac{1+\xi_{K_j,p_j-1}}{\sqrt{1+\frac13\xi_{K_j,p_j-1}^2}+\sqrt2\xi_{K_j,p_j-1}}\in\left(\frac{\sqrt{3}}{\sqrt6+1},1\right),
\end{equation}
cf.~\eqref{eq:range}, to estimate the smoothness of~$f|_{K_j}$. Here, we emphasise that the computation of~$\xi_{K_j,p_j-1}$ does not require any additional function evaluations of~$f$ since the values~$(f\circ\phi_{K_j})(\x{p_j}{k})$, $k=1,\ldots,p_j$, have already been determined in the application of the quadrature rule~\eqref{eq:QKj}.

\subsection{Adaptive Variable Order Procedure}\label{sc:hprefine}

Based on the above derivations, we now propose an $hp$-type adaptive quadrature method. To this end, we start by choosing a tolerance~$\tol>0$, a smoothness parameter~$\tau\in\left(\nicefrac{\sqrt{3}}{(\sqrt6+1)},1\right)$, and a maximal number~$p_{\max}\ge 2$ of possible quadrature points on each subinterval. Furthermore, we define the interval $K_1=[a,b]$, and a small number~$p_1$, $2\le p_1\le p_{\max}$, of quadrature points on~$K_1$. Moreover, we initialise the set of subintervals~$\subs$, the order vector~$\p$ containing the number of quadrature points on each subinterval, and the unknown value~$\Q$ of the integral as follows:
\[
\subs=\{K_1\},\qquad \p=\{p_1\},\qquad \Q=0.
\]
Then, the basic adaptive procedure is given as follows:
\begin{algorithmic}[1]
\While {$\subs\neq\emptyset$}
\State $[\Q1,\subs,\p] = \mathtt{hprefine}(f,\subs,\p,p_{\max},\tau)$;
\State $\Q = \Q + \Q1$;
\EndWhile 
\State Output~$\Q$.
\end{algorithmic}
Here, $\mathtt{hprefine}$ is a function, whose purpose is to identify those subintervals in~$\subs$, which need to be refined further for a sufficiently accurate approximation of the unknown integral. In addition, it outputs a set of subintervals (again denoted by~$\subs$), as well as an associated order vector (again denoted by~$\p$) which result from applying the most appropriate refinement, i.e., either~$h$- or $p$-refinement as outlined in~(i) and~(ii) in Section~\ref{sc:basic} above, for each subinterval. Furthermore, $\mathtt{hprefine}$ returns the sum~$\Q1$ of all quadrature values corresponding to subintervals in the input set~$\subs$ for which no further refinement is deemed necessary. The essential steps are summarised in Algorithm~\ref{alg:hprefine}.

\begin{algorithm}
\caption{Function $[\Q,\subsnew,\pnew] = \mathtt{hprefine}(f,\subs,\p,p_{\max},\tau)$}
\label{alg:hprefine}
\begin{algorithmic}[1]
\State {Define~$\subsnew=\subs$, and~$\pnew=\p$. Set~$\Q=0$.}
\For {each subinterval~$K_j\in\subs$}
\State {Evaluate the smoothness indicator~$\F_{K_j,p_j}(f)$ from~\eqref{eq:sind}.}
\If {$\F_{K_j,p_j}(f)<\tau$}
\myStateDouble {Apply $h$-refinement to~$K_j$, i.e., bisect~$K_j$ into two subintervals of equal size and reduce the number of quadrature points to~$\max(p_j-1,1)$ on both of them;}
\myStateDouble {Compute an improved approximation, denoted by~$\widetilde Q_{K_j}$, of~$Q_{K_j,p_j}(f|_{K_j})$ using~\eqref{eq:Qh} on~$K_j$.}
\ElsIf {$\F_{K_j,p_j}(f)\ge\tau$ and $p_j+1\le p_{\max}$}
\myStateDouble {Apply $p$-refinement to~$K_j$, i.e., increase the number of quadrature points to~$p_j+1$ on~$K_j$;}
\myStateDouble {Compute an improved approximation, denoted by~$\widetilde Q_{K_j}$, of~$Q_{K_j,p_j}(f|_{K_j})$ using~\eqref{eq:Qp} on~$K_j$.}
\ElsIf {$\F_{K_j,p_j}(f)\ge\tau$ and $p_j+1>p_{\max}$}
\myStateDouble {Bisect~$K_j$ into two subintervals of equal size and retain the number of quadrature points~$p_j$ on both of them;}
\myStateDouble {Compute an improved approximation, denoted by~$\widetilde Q_{K_j}$, of~$Q_{K_j,p_j}(f|_{K_j})$ using~\eqref{eq:Qp'} on~$K_j$.}
\EndIf
\If {$|\widetilde Q_{K_j}-Q_{K_j,p_j}(f|_{K_j})|$ is sufficiently small}
\myStateDouble {Update $\Q = \Q + \widetilde Q_{K_j}$;}
\myStateDouble {Eliminate $K_j$ from~$\subsnew$ and the corresponding entry~$p_j$ from~$\pnew$.}
\Else 
\myStateDouble {Replace $K_j$ and~$p_j$ in~$\subsnew$ and~$\pnew$, respectively, by the corresponding $h$- or $p$-refined subintervals as determined above.}
\EndIf 
\EndFor
\end{algorithmic}
\end{algorithm}

\subsection{Practical Aspects}
In this section we discuss a number of practical issues involved in the implementation of the procedure described in Section~\ref{sc:hprefine} within a given computing environment.

\subsubsection{Gauss-Quadrature Rules}
In principle, the adaptive procedure presented in Section~\ref{sc:hprefine} allows for any variable order family of quadrature rules. In our numerical experiments presented in Section~\ref{sc:numerics} below, we propose the use of (families of) Gauss-type quadrature schemes. Although they might be criticised for their non-hierarchical structure, in the sense that they require more function evaluations in comparison to more traditional schemes (such as, for example, the adaptive Simpson or fixed-order Gauss-Kronrod rules), our numerical results indicate that their high degree of accuracy may be exploited in a very efficient manner within the $hp$-setting, particularly for smooth functions, with or without locally singular behaviour. Indeed, whilst non-hierarchical lower-order Gauss-type quadrature schemes might not be computationally competitive, it is a well-known feature of $hp$-methods (see, e.g., \cite{Schwab98}) that their superiority becomes especially apparent on a variable, higher-order level. 

In the current article we employ Gauss-Legendre quadrature points and weights (with at least~$p_{\min}=2$ points and weights); these quantities can be precomputed up to any given order~$p_{\max}$ (in practice~$p_{\max}=15$ is usually more than sufficient) or even be generated on the spot in an efficient way (see, e.g.,~\cite{CanutoHussainiQuarteroni:88,GlaserLiuRokhlin:07,Waldvogel:06}) if an upper bound~$p_{\max}$ cannot be fixed. In addition, we note that the Gauss-Legendre rule based on~$p$ points has a degree of exactness of~$2p-1$, i.e., the smoothness indicators derived in Section~\ref{sc:smoothness} can be computed by means of the formula given in~\eqref{eq:xi}. For a given maximum number~$p_{\max}$, we store the points and weights of the Gauss-Legendre rules (on the reference interval~$[-1,1]$) with up to~$p_{\max}$ points in two $p_{\max}\times (p_{\max}-1)$-matrices~$\bm X$ and~$\bm W$, respectively; here, for parameters~$p=2,\ldots,p_{\max}$, the $p$-th columns of~$\bm X$ and~$\bm W$ are built from the points and weights of the corresponding $p$-point Gauss-Legendre quadrature rule, respectively (and complementing the remaining entries in all but the last column by zeros):
\begin{equation}\label{eq:XW}
{\bm X}=
\begin{pmatrix}
    \x{2}{1}&\x{3}{1}& \cdots & \x{p_{\max}}{1} \\
    \x{2}{2}&\vdots  &        &  \\
            &\x{3}{3}&        & \vdots\\
            &\bigzero& \ddots & \\
            &        &        & \x{p_{\max}}{p_{\max}}
\end{pmatrix},\quad
{\bm W}=
  \begin{pmatrix}
    \w{2}{1}&\w{3}{1}& \cdots & \w{p_{\max}}{1} \\
    \w{2}{2}&\vdots  &        &  \\
            &\w{3}{3}&        & \vdots\\
            &\bigzero& \ddots & \\
            &        &        & \w{p_{\max}}{p_{\max}}
  \end{pmatrix}.
\end{equation}
We note that, for other quadrature rules, the number of rows in the above matrices may be different.

\subsubsection{Vectorised Quadrature}\label{sc:vecQuad}
Following the ideas of~\cite{Shampine:08} we use a vectorised quadrature implementation. This means that, instead of computing the integrals on the subintervals~$\subs$ in Algorithm~\ref{alg:hprefine} one at a time, they are all computed at once. This can be accomplished by using fast vector- and matrix-operations, and by carrying out all necessary function evaluations in a single operation by computing the function to be integrated for a vector of input values. Specifically, we write the composite rule
\[
I\approx \sum_{K_i\in\subs}Q_{K_i,p_i}(f|_{K_i})=\sum_{K_i\in\subs}\frac{h_i}{2}\sum_{k=1}^{p_i}\w{p_i}{k}(f\circ\phi_{K_i})(\x{p_i}{k})
\]
as a dot product of a weight vector~$\bm w$ and a function vector~$f(\bm x)$; here, the former vector contains all (scaled) weights~$\{\frac12h_i\w{p_i}{k}\}_{i,k}$, and the latter vector represents the evaluation of the integrand function~$f$ on the vector~$\bm x$ of all corresponding quadrature points~$\{\phi_{K_i}(\x{p_i}{k})\}_{i,k}$ appearing in the sum above. Evidently, these vectors can be built efficiently by extracting (and affinely mapping and scaling) the corresponding rows from the matrices~$\bm X$ and~$\bm W$ in~\eqref{eq:XW}. We emphasise that applying vectorised quadrature crucially improves the performance of the overall adaptive procedure (provided that such a technology is available in a given computing environment). 

\subsubsection{Smoothness Estimators}
As mentioned before, computing the smoothness indicators from \eqref{eq:xi} does not need any additional function evaluations of the integrand function~$f$; they only require the values of the Legendre polynomials~$\widehat L_{p-1}$ and~$\widehat L_{p-2}$ at the points~$\{\x{p}{k}\}_{k=1}^p$, for~$p=2,\ldots,p_{\max}$. These quantities are again precomputable, and can be stored in two matrices
\begin{equation}\label{eq:L1}
\bm L_1
=\begin{pmatrix}
L_1(\x{2}{1}) & L_2(\x{3}{1}) & \cdots & L_{p_{\max}-1}(\x{p_{\max}}{1})\\
L_1(\x{2}{2}) & \vdots\\
              & L_2(\x{3}{3}) & & \vdots\\
              &\bigzero& \ddots & \\
              &        && L_{p_{\max}-1}(\x{p_{\max}}{p_{\max}})
\end{pmatrix},
\end{equation}
and
\begin{equation}\label{eq:L2}
\bm L_2
=\begin{pmatrix}
L_0(\x{2}{1}) & L_1(\x{3}{1}) & \cdots & L_{p_{\max}-2}(\x{p_{\max}}{1})\\
L_0(\x{2}{2}) & \vdots\\
              & L_1(\x{3}{3}) & & \vdots\\
              &\bigzero& \ddots & \\
              &        && L_{p_{\max}-2}(\x{p_{\max}}{p_{\max}})
\end{pmatrix}.
\end{equation}
Then, the sums in~\eqref{eq:xi} are vectorised similarly as described above. In particular, the computation of the smoothness estimators can be undertaken with an almost negligible computational cost.

\subsubsection{Stopping Criterion}
In order to implement the stopping-type criterion in line~14 of Algorithm~\ref{alg:hprefine}, we exploit an idea that was proposed in the context of adaptive Simpson quadrature in~\cite{GanderGautschi:00}. More precisely, given a possibly rough approximation~$\mathtt{iguess}\approx\int_a^bf(x)\dx$ of the unknown integral~$I$ from~\eqref{eq:I} (e.g., obtained from a Monte-Carlo calculation such that both the approximation and the exact value are of the same magnitude; cf.~\cite{GanderGautschi:00}), and a tolerance~$\tol>0$, we redefine
\[
\mathtt{iguess = iguess*tol/eps;}
\]
here, $\mathtt{eps}$ represents the smallest (positive) machine number in a given computing environment. Then, using the comparison operator~$\mathtt{==}$, we accept the difference $|\widetilde Q_{K_j}-Q_{K_j,p_j}(f|_{K_j})|$ to be sufficiently small with respect to the given tolerance~$\tol$ if the logical call
\[
\mathtt{iguess} + |\widetilde Q_{K_j}-Q_{K_j,p_j}(f|_{K_j})| \mathtt{\ == iguess};
\]
yields a {\tt true} value. 

\subsection{Numerical Examples}\label{sc:numerics}
\begin{figure}
\centering
\includegraphics[width=0.685\linewidth]{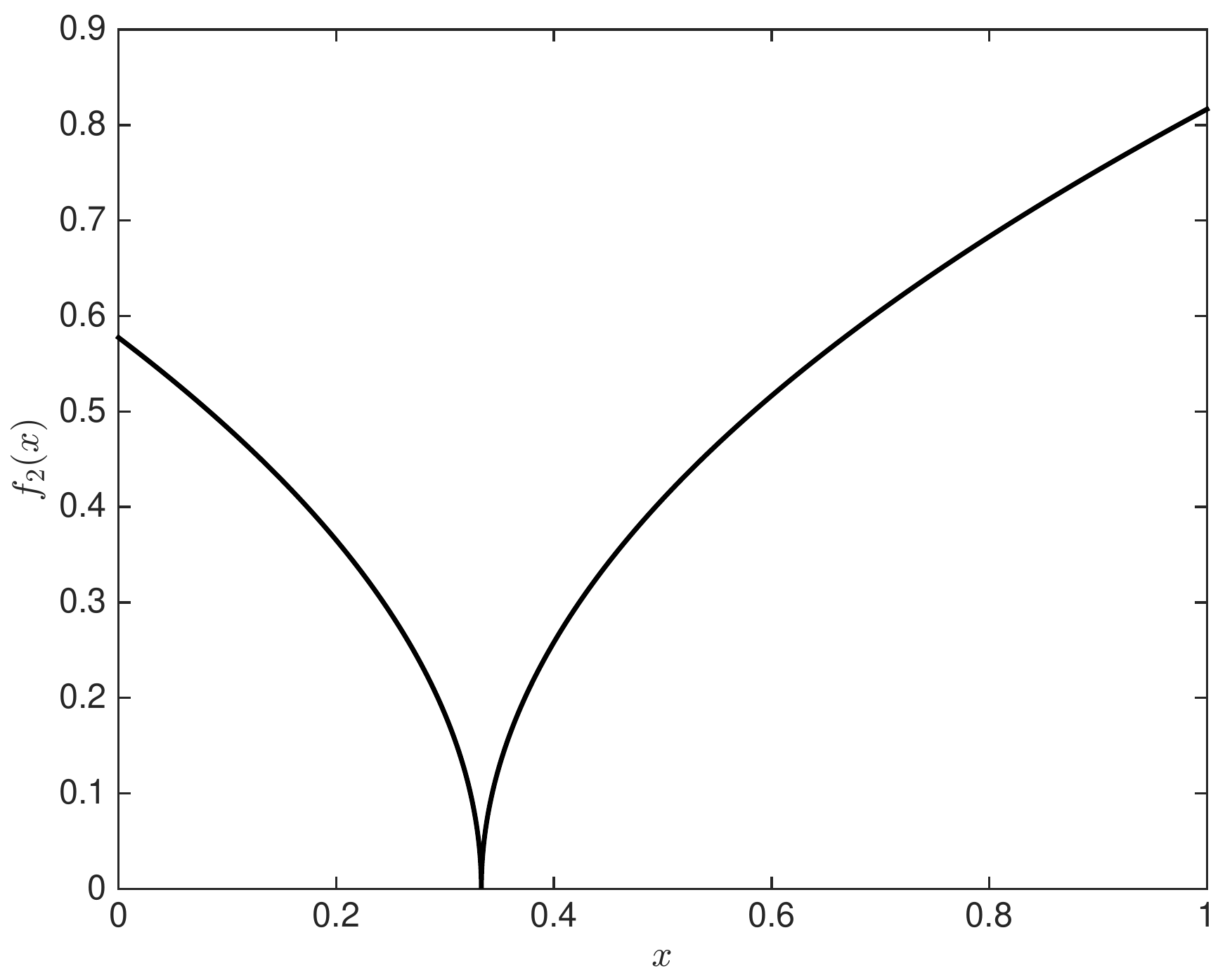}\\[5ex]
\includegraphics[width=0.685\linewidth]{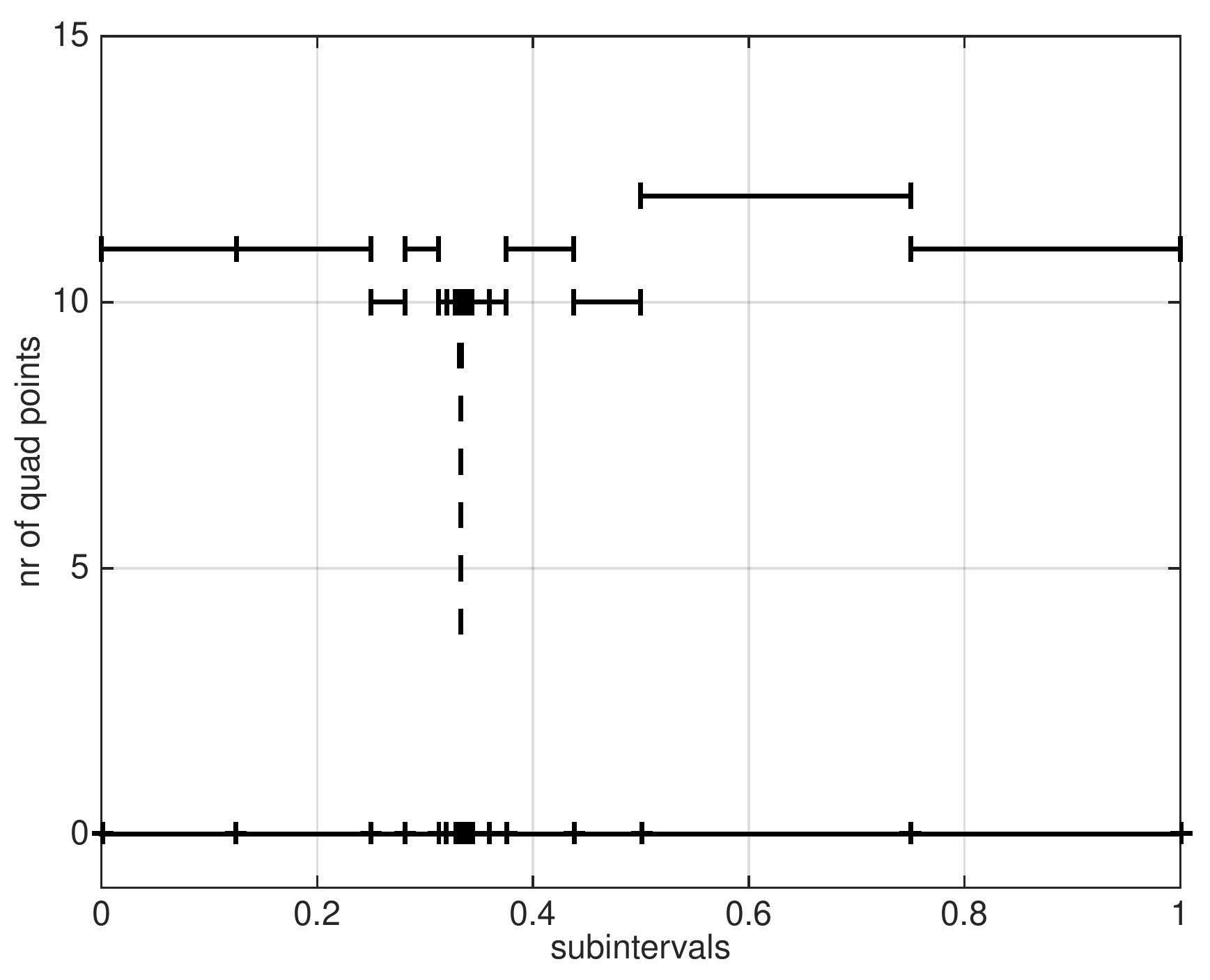}
\caption{Function~$f_2$: Graph (top) and $hp$-mesh (bottom).}
\label{fig:f2}
\end{figure}

\begin{figure}
\centering
\includegraphics[width=0.685\linewidth]{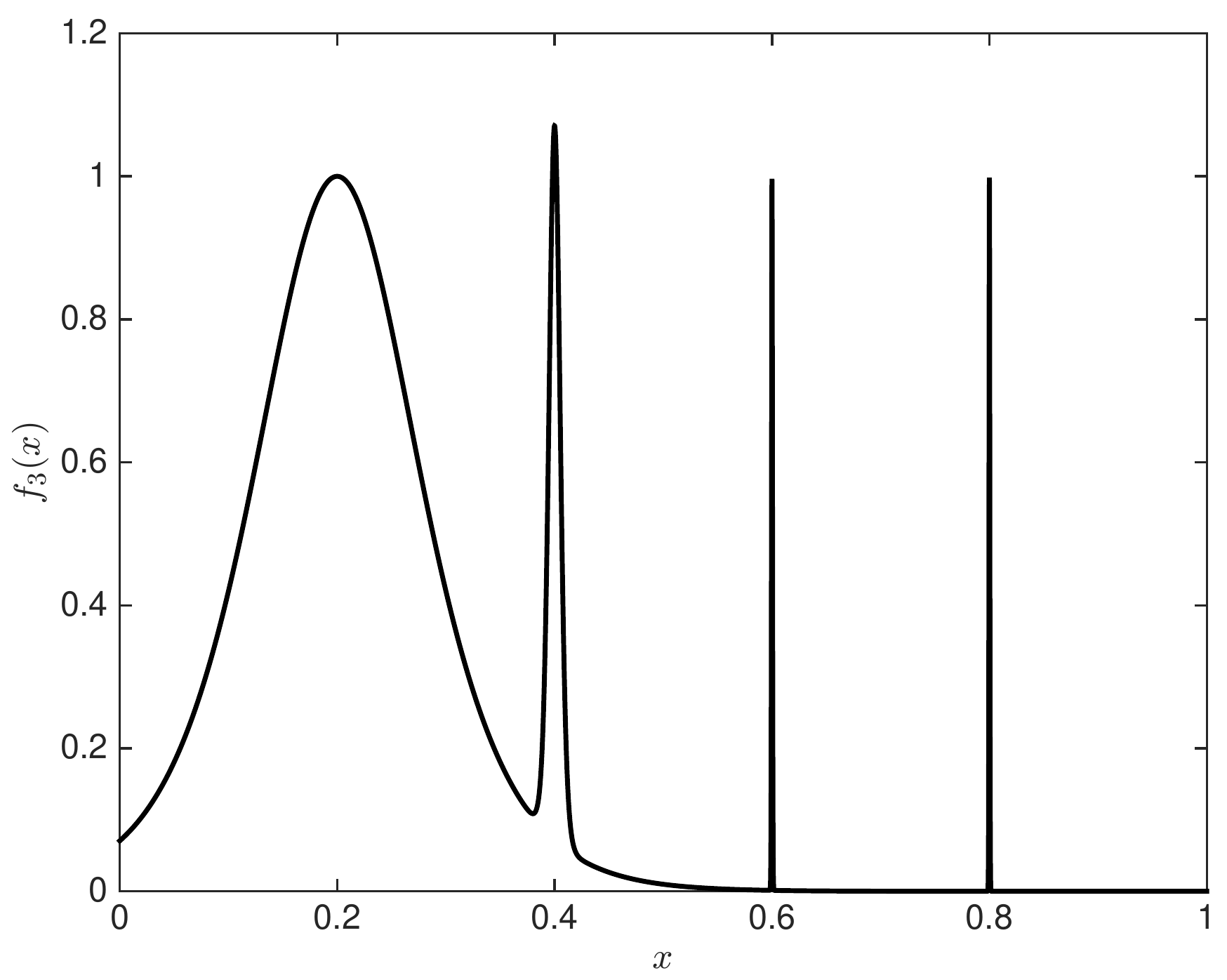}\\[5ex]
\includegraphics[width=0.685\linewidth]{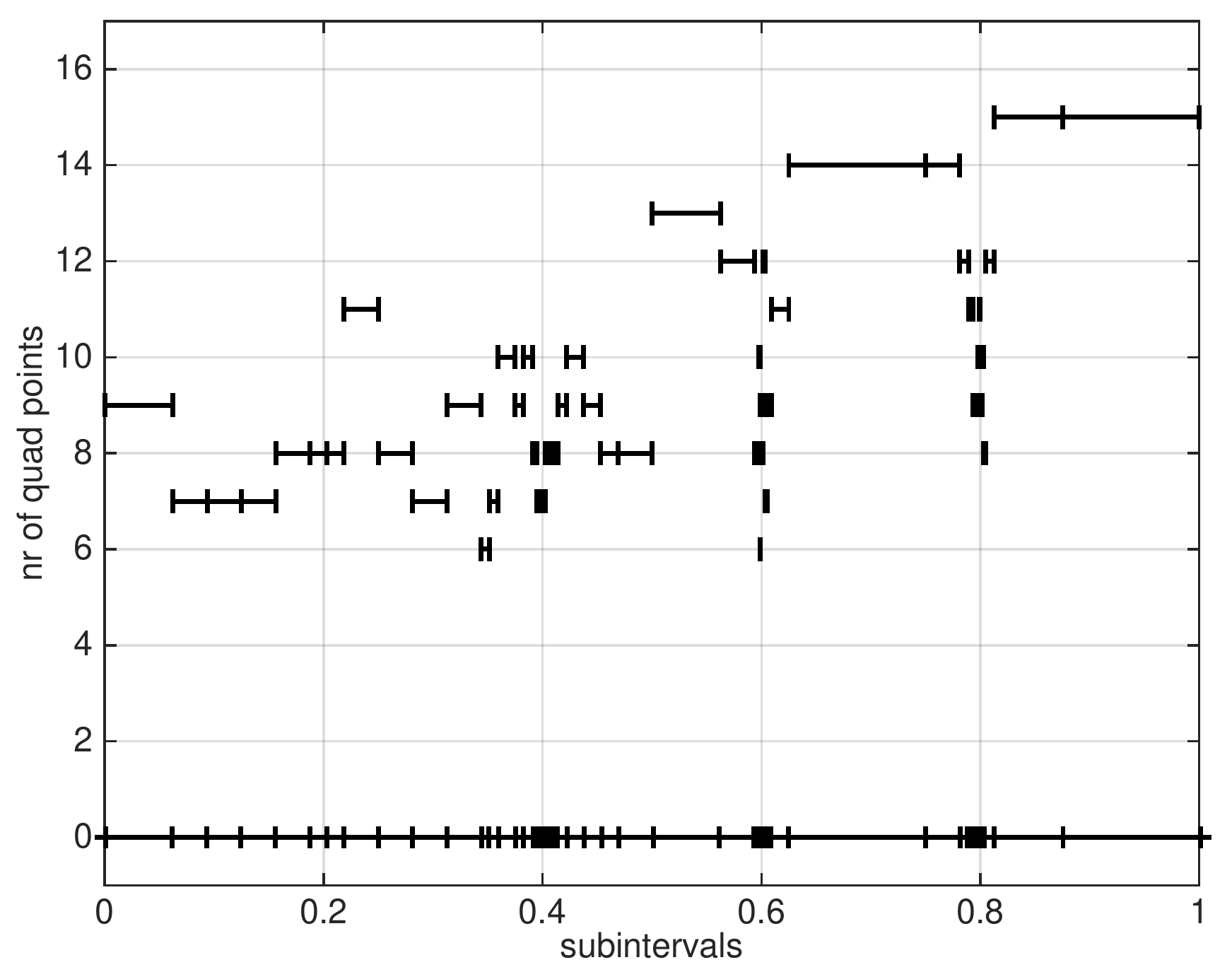}
\caption{Function~$f_3$: Graph (top) and $hp$-mesh (bottom).}
\label{fig:f3}
\end{figure}

In order to test our approach, we consider a number of benchmark problems on the interval~$[0,1]$. Specifically, the following functions will be studied:
\begin{align*}
f_1(x) & = \exp(x),\\
f_2(x) & = \sqrt{|x-\nicefrac13|},\\
f_3(x) & = \sech(10(x-\nicefrac15))^2 + \sech(100(x-\nicefrac25))^4\\
	 &\quad + \sech(1000(x-\nicefrac35))^6 + \sech(1000(x-\nicefrac45))^8,\\
f_4(x) & = \cos(1000x),\\
f_5(x) & = \begin{cases}
0 & \text{ if $x\le\nicefrac13$},\\
1 & \text{ if $x>\nicefrac13$}.
\end{cases}
\end{align*}
Whilst the first function, $f_1$, is analytic, the second function, $f_2$, is smooth except at~$\nicefrac13$ (see Figure~\ref{fig:f2} (top)). Furthermore, $f_3$ was proposed in~\cite{Hale:10} in the context of the chebfun package~\cite{HaleTrefethen:12}; this is a smooth function that exhibits several very thin spikes (see Figure~\ref{fig:f3} (top)). Moreover, $f_4$ is highly oscillating, and~$f_5$ is an example of a discontinuous function. 

We perform our computations in \matlab\footnote{The MathWorks, Inc.} on a single 2.6GHz processor. The tolerance is set to~$\mathtt{tol}=0.3\times10^{-15}$ (which is close to machine precision in \matlab), the smoothness estimation parameter is prescribed as~$\tau=0.6$, and~$p_{\max}=15$. Within this setting, the adaptive procedure generates results that are accurate to machine precision, for all of the considered examples. In Table~\ref{tb:fctevl}, for each of the functions~$f_1,\ldots,f_5$ above, we present the number of function calls  (\#~fct.~calls) in the vectorised quadrature implementation (counting a single application of the integrand function to a vector input as 1; cf.~Section~\ref{sc:vecQuad}), as well as the number of single function evaluations (\#~sing.~fct.~ev.) taking into account the number of scalar entries of a vector input in each function call. The latter number is compared with the number of scalar function evaluations performed in a classical adaptive Simpson procedure as proposed in~\cite{GanderGautschi:00} (which is based on employing the two end points as well as the midpoint on each subinterval, and reuses the former two points without recomputing). Except for the last function, $f_5$, where a low-order quadrature rule is more effective, the remarkable efficiency of the proposed $hp$-type quadrature becomes clearly visible. This is confirmed with the expeditious cpu times (which do not include the computation of the precomputable matrices~$\bm X, \bm W, \bm L_1,\bm L_2$ from~\eqref{eq:XW}, \eqref{eq:L1}, and \eqref{eq:L2}) for each of the examples.

\begin{table}
\begin{tabular}{crrrr}    \toprule
	\sp	     & \multicolumn{3}{c}{\emph{$hp$-adapt. quad.}} & \emph{adapt. Simpson quad.}\\
		        & \#~fct.~calls & \# sing.~fct.~ev. & cpu [sec] & \# sing.~fct.~ev.\\\midrule
$f_1$    &     52 &   9 & 0.0031 &     4,096\\
$f_2$    &  1,718 &  65 & 0.0224 &    25,488\\
$f_3$    &  2,427 &  33 & 0.0144 &    72,528\\
$f_4$    & 50,534 &  35 & 0.0180 & 1,965,376\\
$f_5$    &  1,273 & 106 & 0.0342 &       784\\
\bottomrule\\
\end{tabular}
\caption{Performance data for $hp$-type adaptive quadrature.}
\label{tb:fctevl}
\end{table}

In order to illustrate how the $hp$-adaptive procedure performs, we depict the final $hp$-mesh for~$f_2$ and~$f_3$ in Figure~\ref{fig:f2} (bottom) and Figure~\ref{fig:f3} (bottom), respectively. Here, along the horizontal axis we present the subintervals obtained as a result of the adaptive process, and on the vertical axis the number of quadrature points introduced on each subinterval is displayed. In both examples, we see that smooth regions in the underlying integrand are resolved by employing larger subintervals featuring a higher number of quadrature points, whereas close to singularities, the number of quadrature points is kept low on very small integration subdomains. It is noteworthy that this behaviour is well-known from $hp$-finite element methods for differential equations, where high-order algebraic or even exponential convergence rates can be obtained by applying this type of $hp$-refinement procedure; see~\cite{Schwab98} for details.

\section{Conclusions}
In this article we proposed a new adaptive quadrature strategy, which features both local subdivision of the integration domain, as well as local variation of the number of quadrature points employed on each subinterval. Our approach is inspired by the $hp$-adaptive finite element methodology based on $hp$-adaptive smoothness testing. In combination with a vectorised quadrature implementation, the proposed adaptive quadrature algorithm is able to deliver highly accurate results in a very efficient manner. Since our approach is closely related to the $hp$-finite element technique, it can be extended to multiple dimensions, including, in particular, the application of anisotropic refinements of the underlying domain of integration, together with the exploitation of different numbers of quadrature points in each coordinate direction on each subinterval (based, for example, on anisotropic Sobolev embeddings as outlined in~\cite[\S3.1]{FankhauserWihlerWirz:14}).

\bibliographystyle{amsplain}
\bibliography{literature}

\end{document}